\documentclass[a4paper,10pt]{amsart}

\usepackage{amsmath,amssymb,amsthm,a4wide}

\usepackage{tikz}
\usepackage{graphicx, caption}
\usepackage{graphics}

\newtheorem*{thm}{Theorem}

\newcommand{\sgn}{\operatorname{sgn}}

\begin{document}

\title[]{An amusing sequence of functions}

\author[]{Stefan Steinerberger}
\address{Department of Mathematics, Yale University, 06511 New Haven, CT, USA}
\email{stefan.steinerberger@yale.edu}

\begin{abstract}  We consider the amusing sequence of functions $f_n: \mathbb{R} \rightarrow \mathbb{R}$ given by
$$ f_n(x) = \sum_{k=1}^{n}{\frac{|\sin{(k \pi x)}|}{k}}.$$
Every rational point is eventually the location of a strict local minimum of $f_n$: more precisely,
$f_n$ has a strict local minimum in all rational points $x=p/q \in \mathbb{Q}$ with $|q| \leq \sqrt{n}$.
\end{abstract}

\maketitle
\vspace{-10pt}
\section{Introduction}

The purpose of this short note is to introduce 
$$ f_n(x) =  \sum_{k=1}^{n}{\frac{|\sin{(k \pi x)}|}{k}}.$$

\begin{thm} The function $f_n(x)$ has a strict local minimum in $x=p/q$ for all $n \geq q^2$.
\end{thm}

The asymptotically sharp scaling is given by $n \geq (1+o(1))q^2/\pi$.
We believe that this curious result is a good indicator that this sequence might have all sorts of other nice
properties and could be of some interest. A natural question would be whether anything can
be said about the location of local maxima: it is tempting to conjecture that they cannot be too well 
approximated by rationals with small denominators
(because that's where the minima are).

\begin{figure}[h!]
\begin{minipage}{0.49\textwidth}
\begin{center}
\begin{tikzpicture}
\node at (0,0) {\includegraphics[width= 0.8\textwidth]{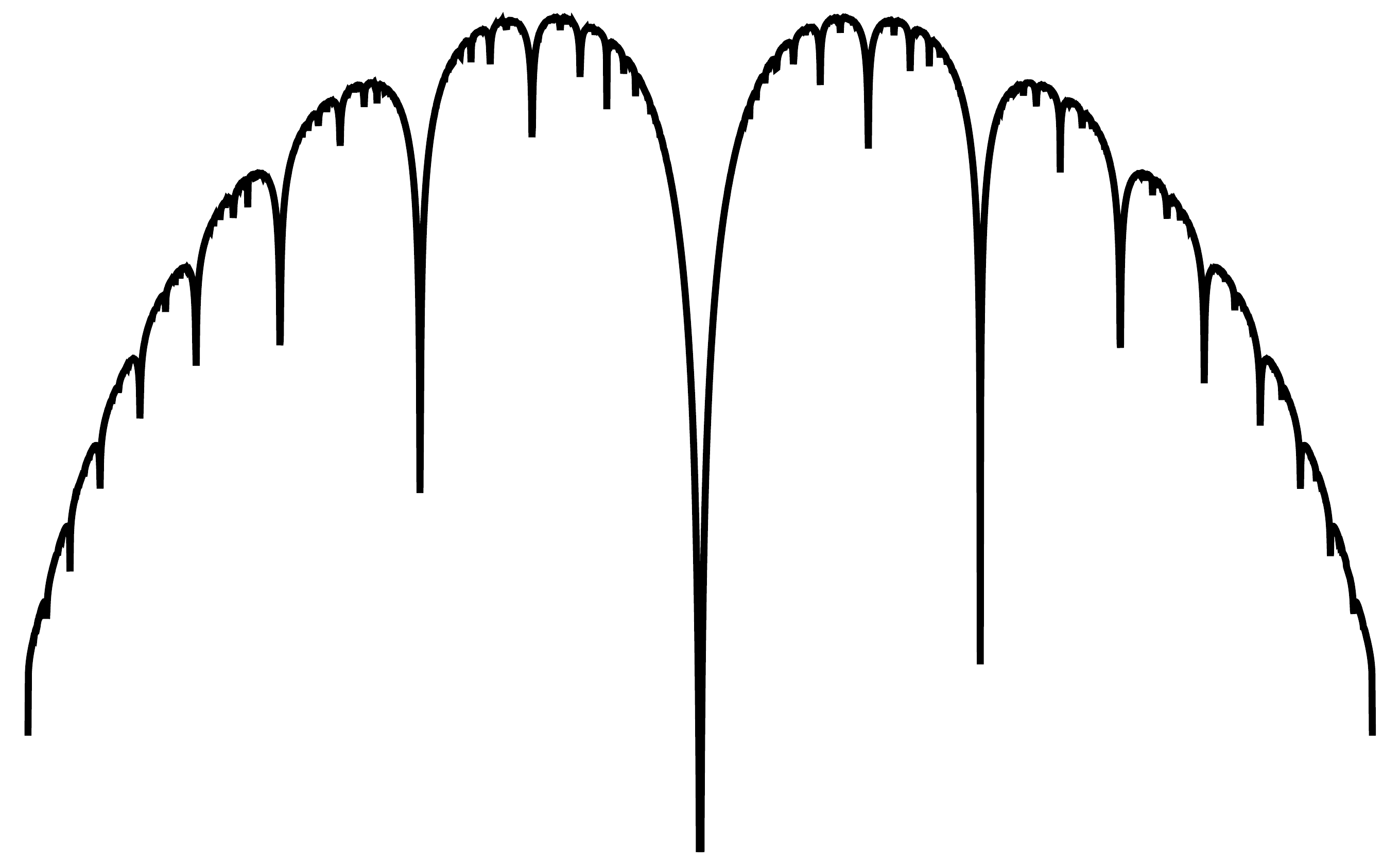}};
\draw [ultra thick] (-3,-2) -- (3,-2);
\draw [ultra thick] (-3,-2.1) -- (-3,-1.9);
\node at (-3, -2.4) {$0.1$};
\draw [ultra thick] (3,-2.1) -- (3,-1.9);
\node at (3, -2.4) {$0.9$};
\end{tikzpicture}
\end{center}
\end{minipage}
\begin{minipage}{0.49\textwidth}
\begin{center}
\begin{tikzpicture}
\node at (0,0) {\includegraphics[width= 0.8\textwidth]{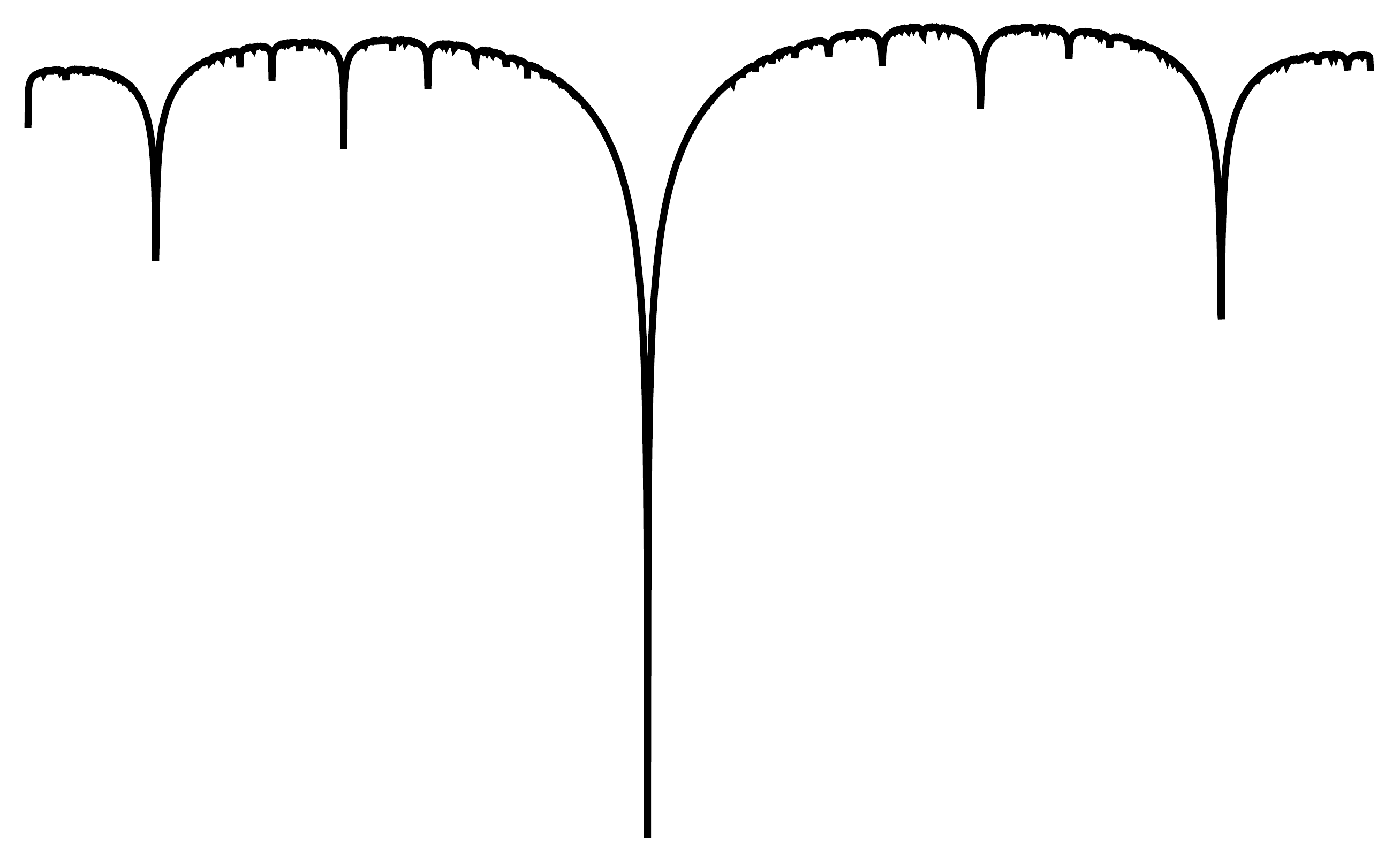}};
\draw [ultra thick] (-3,-2) -- (3,-2);
\draw [ultra thick] (-3,-2.1) -- (-3,-1.9);
\node at (-3, -2.4) {$0.38$};
\draw [ultra thick] (3,-2.1) -- (3,-1.9);
\node at (3, -2.4) {$0.39$};
\end{tikzpicture}
\end{center}
\end{minipage}
\captionsetup{width=0.9\textwidth}
\caption{The function $f_{50.000}$ on $[0.1, 0.9]$ and zoomed in (right). The big cusp in the right picture is located at $x=5/13$,
the two smaller cusps are at $x=8/21$ and $x=7/18$.}
\end{figure}
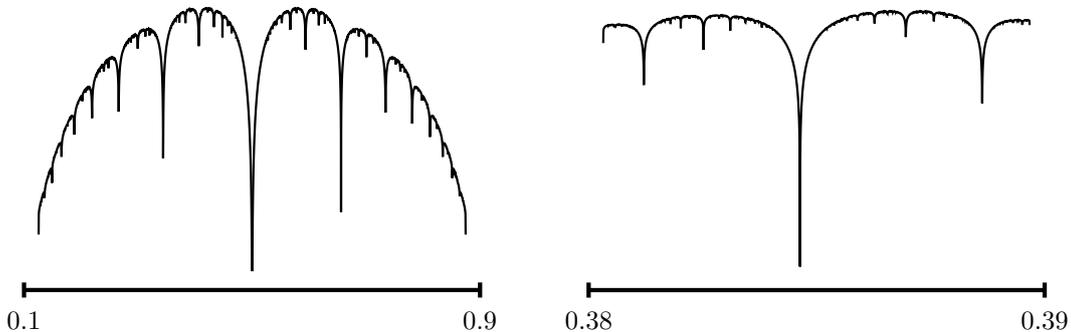

Both definition of the function as well as its graph are reminiscent of the Takagi function $\tau$, which is a continuous but nowhere differentiable function that was first considered by Takagi  \cite{tak} in 1901 (with independent re-discoveries by van der Waerden \cite{vdw} in 1930 and de Rham \cite{der} in 1957): if $d(x)$ denotes the distance from $x$ to the nearest integer, then $\tau$ is given by
$$ \tau(x) = \sum_{k=0}^{\infty}{\frac{d(2^k x)}{2^k}}.$$
$\tau(x)$ has since appeared in connection to inequalities for digit sums \cite{allard}, the Riemann hypothesis \cite{bal} and extremal combinatorics \cite{frankl} (many more results can be found in the surveys \cite{all, lag}). We emphasize a 1959 result of Kahane \cite{kahane} who proved that the set of local minima are
exactly the dyadic rational numbers.

\begin{figure}[h!]
\begin{minipage}{0.49\textwidth}
\begin{center}
\begin{tikzpicture}
\node at (0,0) {\includegraphics[width= 0.8\textwidth]{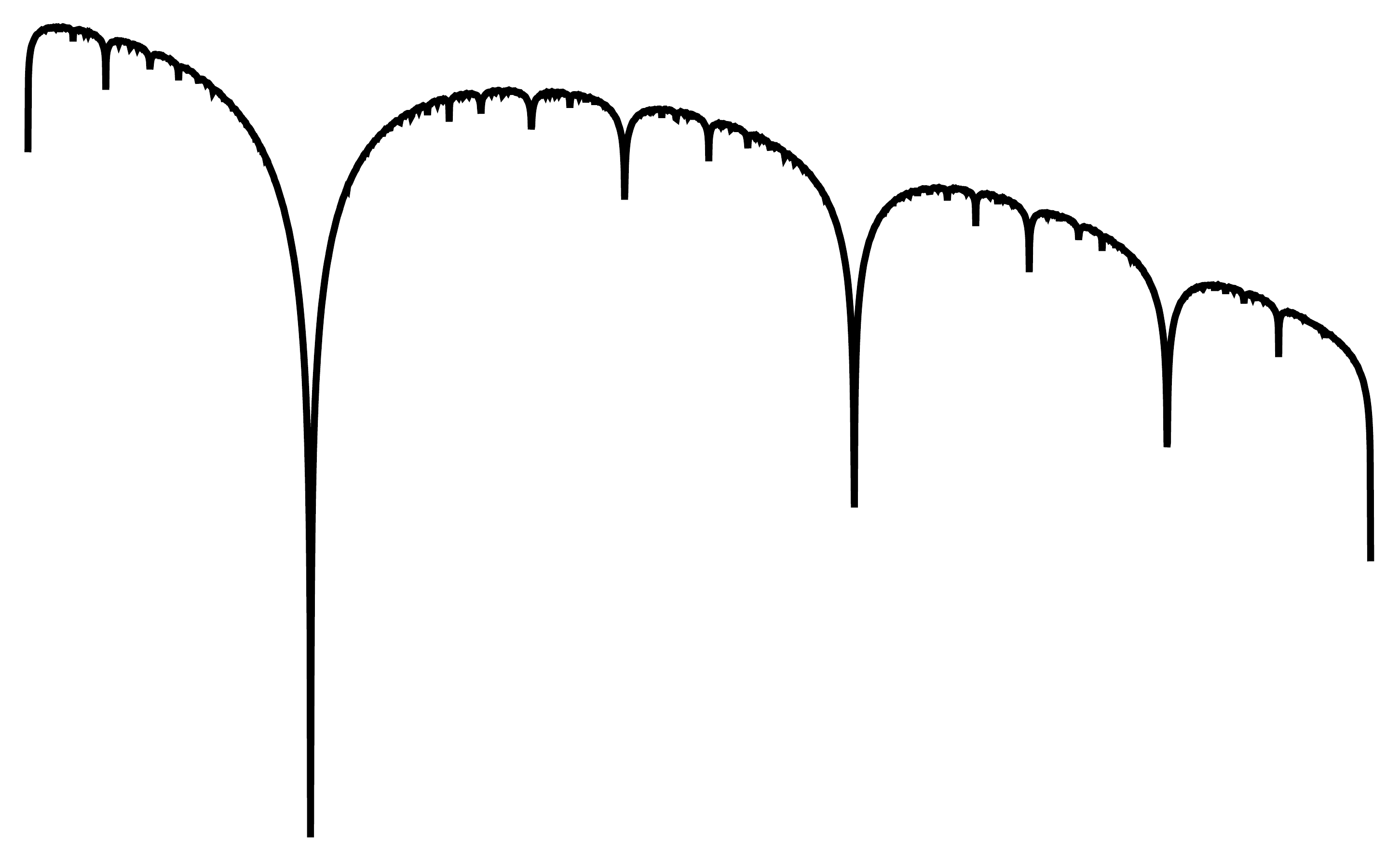}};
\draw [ultra thick] (-3,-2) -- (3,-2);
\draw [ultra thick] (-3,-2.1) -- (-3,-1.9);
\node at (-3, -2.4) {$0.42$};
\draw [ultra thick] (3,-2.1) -- (3,-1.9);
\node at (3, -2.4) {$0.425$};
\end{tikzpicture}
\end{center}
\end{minipage}
\begin{minipage}{0.49\textwidth}
\begin{center}
\begin{tikzpicture}
\node at (0,0) {\includegraphics[width= 0.8\textwidth]{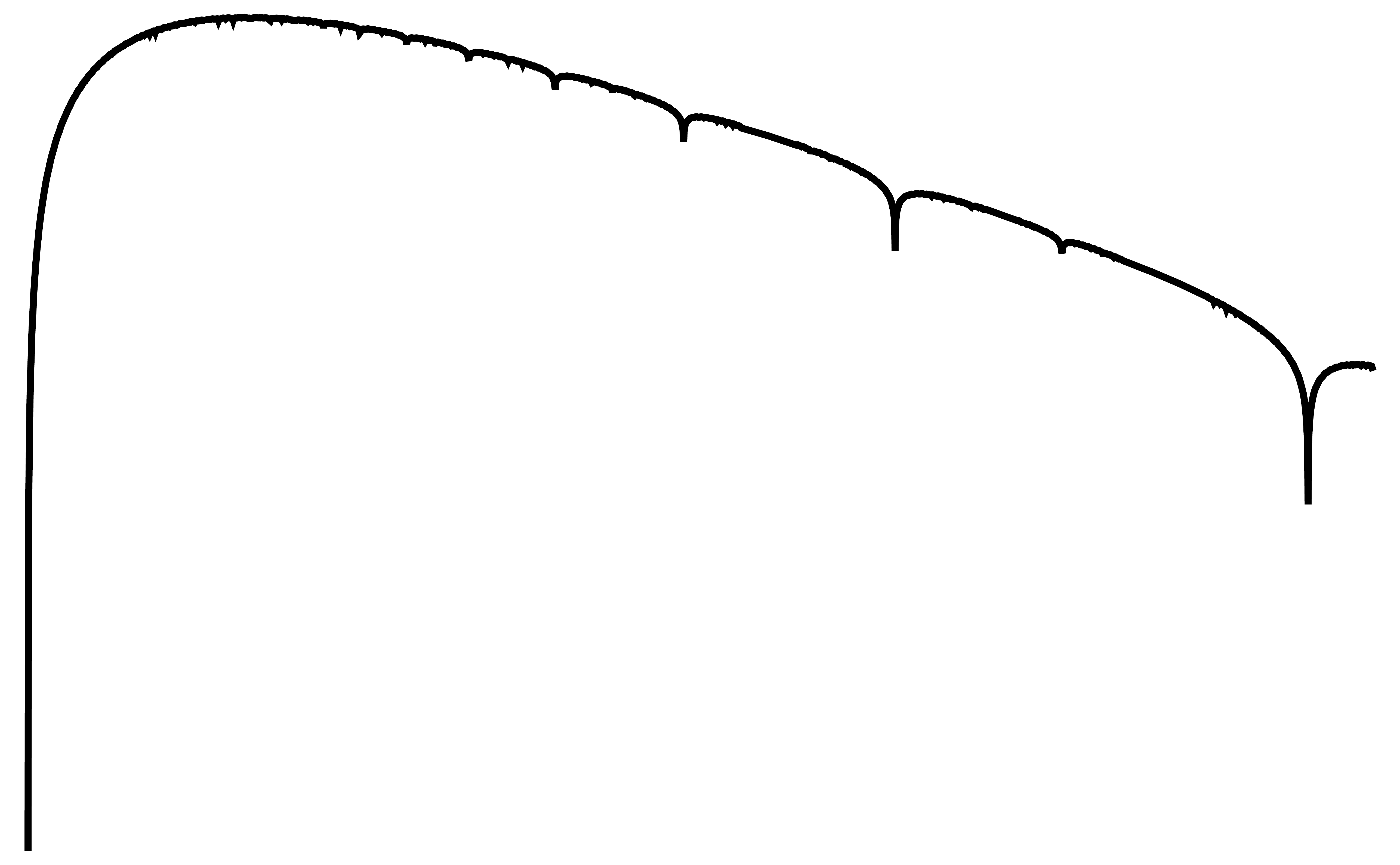}};
\draw [ultra thick] (-3,-2) -- (3,-2);
\draw [ultra thick] (-3,-2.1) -- (-3,-1.9);
\node at (-3, -2.4) {$0.9$};
\draw [ultra thick] (3,-2.1) -- (3,-1.9);
\node at (3, -2.4) {$0.905$};
\end{tikzpicture}
\end{center}
\end{minipage}
\captionsetup{width=0.9\textwidth}
\caption{The function $f_{50.000}$ in two other locations.}
\end{figure}
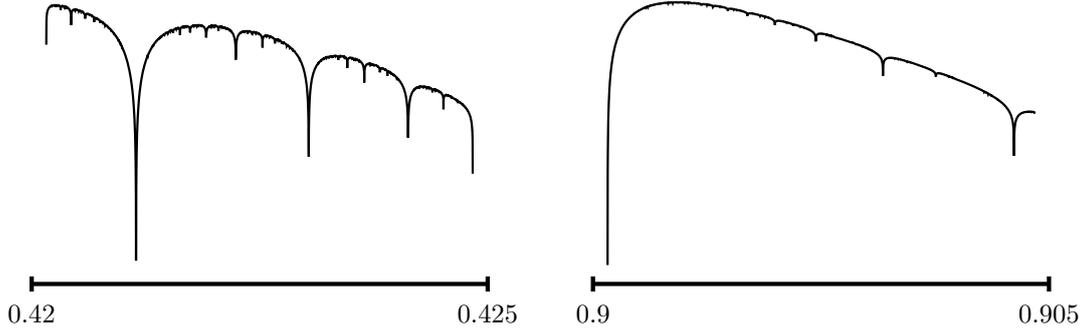

\newpage

\section{Proof} 
\begin{proof} We observe that for $x \in \mathbb{R}$ and $\varepsilon \rightarrow 0$
$$ |\sin{(x+\varepsilon)}| - |\sin{(x)}| = \begin{cases} |\varepsilon| + \mathcal{O}(\varepsilon^2) \qquad &\mbox{if}~x/\pi \in \mathbb{Z} \\
\varepsilon\sgn(\sin{(x)}) \cos{(x)} + \mathcal{O}(\varepsilon^2)  \qquad &\mbox{otherwise.} \end{cases}$$
Here, $\mbox{sgn}$ denotes the signum function
$$ \mbox{sgn}(x) = \begin{cases} 1 \qquad &\mbox{if}~x > 0 \\
0\qquad &\mbox{if}~x = 0\\
-1 \qquad &\mbox{if}~x < 0.
\end{cases}$$
Let us now consider the function in $x = p/q$ with $\mbox{gcd}(p,q) = 1$.
We have
\begin{align*}
\sum_{k=1}^{n}{ \frac{|\sin{(k \pi (p/q + \varepsilon))}|}{k} - \frac{|\sin{(k \pi (p/q))}|}{k} } &=   \pi \varepsilon \sum_{k=1}^{n}{ \sgn\left(\sin{ \left( \frac{k \pi p}{q}\right)}\right) \cos{ \left( \frac{k \pi p}{q}\right)}  } \\
&+  \pi |\varepsilon| \# \left\{1 \leq k \leq n: k (p/q) \in \mathbb{Z} \right\} + \mathcal{O}(\varepsilon^2).
\end{align*}
We analyze these two coefficients and show that the first one is bounded (the second term is clearly unbounded).
The function $\sgn(\sin{(x)}) \cos{(x)}$ has period $\pi$ and the map $k \rightarrow k \cdot p$ is a permutation on $\mathbb{Z}_q$. It is then easy to see that the symmetries of sine and cosine imply
$$\sum_{k=1}^{q}{ \sgn\left(\sin{ \left( \frac{k \pi p}{q}\right)}\right) \cos{ \left( \frac{k \pi p}{q}\right)}  } = 0 \quad \mbox{and thus} \quad \sum_{k=m+1}^{m+q}{ \sgn\left(\sin{ \left( \frac{k \pi p}{q}\right)}\right) \cos{ \left( \frac{k \pi p}{q}\right)}  } = 0$$
for all $m \in \mathbb{N}$.
The periodicity and $k \rightarrow k \cdot p$ being a permutation give
\begin{align*}
\inf_{n \in \mathbb{N}} \sum_{k=1}^{n}{ \sgn\left(\sin{ \left( \frac{k \pi p}{q}\right)}\right) \cos{ \left( \frac{k \pi p}{q}\right)}  } &= \min_{1 \leq n \leq q}{ \sum_{k=1}^{n}{ \sgn\left(\sin{ \left( \frac{k \pi p}{q}\right)}\right) \cos{ \left( \frac{k \pi p}{q}\right)}  }  } \\
&\geq - \max_{1 \leq n \leq q}{ \sum_{k=1}^{n}{ \cos{ \left( \frac{k \pi }{q}\right)}  }  } \geq - \frac{q}{2},
\end{align*}
while, at the same time, we obviously have 
$$  \# \left\{1 \leq k \leq n: k (p/q) \in \mathbb{Z} \right\}  \geq \left\lfloor \frac{n}{q} \right\rfloor$$
from which positivity follows for $n \geq q^2/2$.\\
\vspace{10pt}
\end{proof}
\textit{Remarks.} A more careful analysis shows that
$$ - \max_{1 \leq n \leq q}{ \sum_{k=1}^{n}{ \cos{ \left( \frac{k \pi }{q}\right)}  }  }  = - (1+o(1))  \frac{q}{2}\frac{2}{\pi} \int_{0}^{\pi/2}{\cos{x}dx} = - (1+o(1))\frac{ q}{\pi}$$
from which we get that, asymptotically, $n \geq (1+o(1))q^2/\pi$ suffices. 
We observe that this is optimal: if $p = q-1$ and $1 \leq k \leq q/2$, then
$$  \sgn\left(\sin{ \left( \frac{k \pi p}{q}\right)}\right) \cos{ \left( \frac{k \pi p}{q}\right)} = - \cos{ \left( \frac{k \pi }{q}\right)} \quad \mbox{and} \quad 
\sum_{k=1}^{\lfloor q/2 \rfloor}{ - \cos{ \left( \frac{k \pi }{q}\right)}} \sim -\frac{q}{\pi}.$$

The main argument only appealed to certain fairly elementary symmetry properties of the trigonometric
functions and easily extends to various other functions. A particularly nice example comes from replacing the sine
by the cosine 
$$ g_n(x) = \sum_{k=1}^{n}{\frac{|\cos{(k \pi x)}|}{k}}.$$
Whether $x= p/q$ is the location of a local minimum or maximum now depends on the parity of $q$. What other
phenomena can be found?

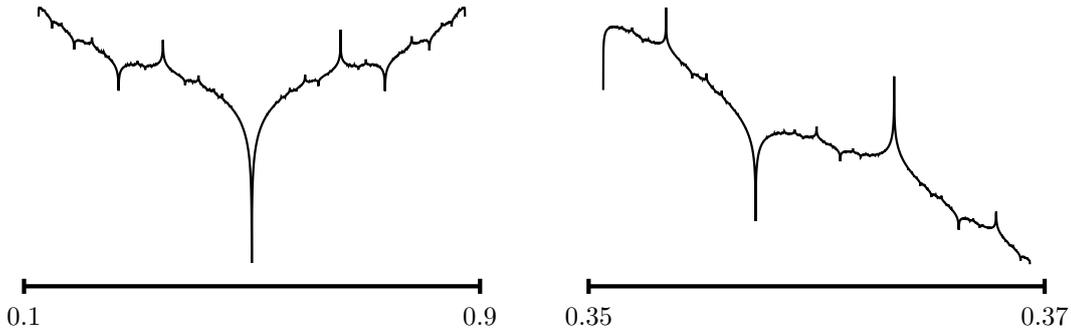
\begin{figure}[h!]
\begin{minipage}{0.49\textwidth}
\begin{center}
\begin{tikzpicture}
\node at (0,0) {\includegraphics[width= 0.8\textwidth]{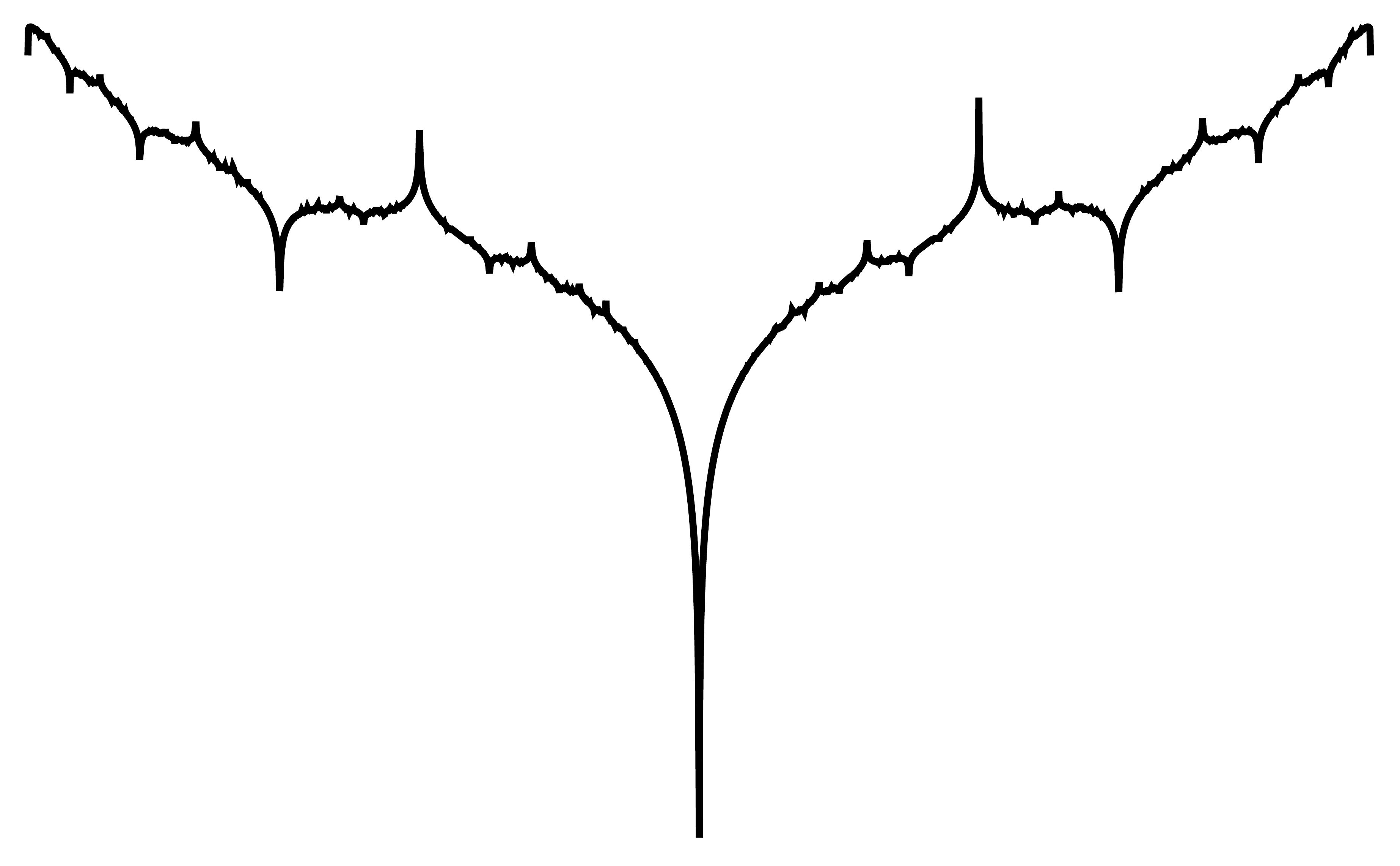}};
\draw [ultra thick] (-3,-2) -- (3,-2);
\draw [ultra thick] (-3,-2.1) -- (-3,-1.9);
\node at (-3, -2.4) {$0.1$};
\draw [ultra thick] (3,-2.1) -- (3,-1.9);
\node at (3, -2.4) {$0.9$};
\end{tikzpicture}
\end{center}
\end{minipage}
\begin{minipage}{0.49\textwidth}
\begin{center}
\begin{tikzpicture}
\node at (0,0) {\includegraphics[width= 0.8\textwidth]{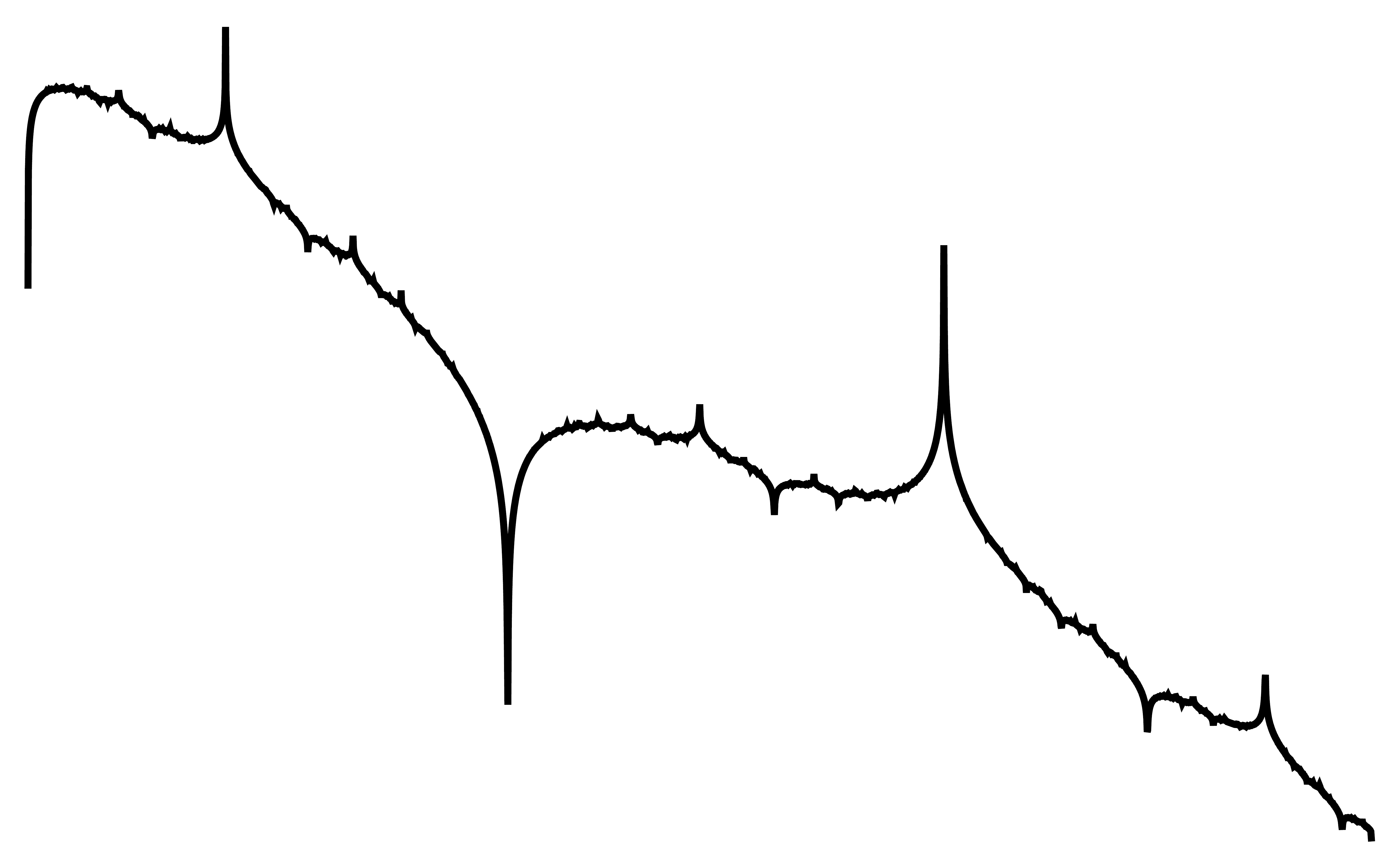}};
\draw [ultra thick] (-3,-2) -- (3,-2);
\draw [ultra thick] (-3,-2.1) -- (-3,-1.9);
\node at (-3, -2.4) {$0.35$};
\draw [ultra thick] (3,-2.1) -- (3,-1.9);
\node at (3, -2.4) {$0.37$};
\end{tikzpicture}
\end{center}
\end{minipage}
\captionsetup{width=0.9\textwidth}
\caption{The function $g_{50.000}$ (left) and zoomed in (right).}
\end{figure}

\textbf{Acknowledgement.} This paper arose out of an entertaining discussion with Raphy Coifman.

\end{document}